\documentclass[a4paper,12pt]{article}

\usepackage{amsmath,amssymb,amsthm}
\usepackage{graphicx,color,hyperref}
\usepackage{rawfonts}
\usepackage{subcaption}

\newcounter{mathitem}
\newenvironment{mathitem}
  {\begin{list}{$(\roman{mathitem})$}{
   \setcounter{mathitem}{0}
   \usecounter{mathitem}
   \setlength{\topsep}{0pt plus 2pt minus 0pt}
   \setlength{\parskip}{0pt plus 2pt minus 0pt}
   \setlength{\partopsep}{0pt plus 2pt minus 0pt}
   \setlength{\parsep}{0pt plus 2pt minus 0pt}
   \setlength{\leftmargin}{35pt}
   \setlength{\itemsep}{0pt plus 2pt minus 0pt}}}
 {\end{list}}

\newtheorem{theorem}{Theorem}

\title{Hamiltonian paths in iterated line graphs}
\author{Jan Ekstein\thanks{Department of Mathematics and European Centre of Excellence NTIS - New Technologies for the Information Society, Faculty of Applied Sciences, University of West Bohemia, Pilsen, Technick\'a 8, 306 14 Plze\v n, Czech Republic, EU
\newline e-mail: \texttt{ekstein@kma.zcu.cz}.}\and
Zuzana Kulhánková \thanks{Department of Mathematics, Faculty of Applied Sciences, University of West Bohemia, Pilsen, Technick\'a 8, 306 14 Plze\v n, Czech Republic, EU
\newline e-mail: \texttt{kulhanzu@kma.zcu.cz}.}}
\date{\today}

\begin{document}
\maketitle

\begin{abstract}
For integer $n$, the $n$-iterated line graph $L^n(G)$ of an undirected graph $G$ is defined to be $L(L^{n-1}(G))$, where $L^1(G)$ is the line graph $L(G)$ of $G$. In this paper we  introduce a hamiltonian path index. The hamiltonian path index, denoted by $h_p(G)$, is the minimum number $n$ such that $L^n(G)$ contains a hamiltonian path. We show that the hamiltonian path index of $G$ exists for any graph $G$ and we set the exact value of the hamiltonian path index for trees and as a generalization for graphs with hamiltonian connected blocks. 

\medskip 
     
\medskip

\end{abstract}

\section{Introduction}
In this paper we consider finite undirected graphs.  As for standard termino\-logy and other terminology used in this paper, we refer to the book by Bondy and Murty, \cite{Bon}. The line graph $L(G)$ of a graph $G$ has edges of $G$ as its vertex set and two vertices are adjacent in $L(G)$ if and only if corresponding edges are adjacent in $G$. For integer $n$, the $n$-iterated line graph $L^n(G)$ of a graph $G$ is defined to be $L(L^{n-1}(G))$, where $L^1(G)=L(G)$, and $L^{n-1}(G)$ is assumed to have a nonempty edge set. We also mean $L^0(G)=G$. We will consider connected graphs only with respect to the study of hamiltonian properties.

Chartrand in \cite{Char} considered the $n$-iterated line graphs and introduced the concept of hamiltonian index of a graph. The hamiltonian index, denoted by $h(G)$, is the minimum number $n$ such that  $L^n(G)$ is hamiltonian. He showed that for any graph $G$ other than a path, the hamiltonian index of $G$ exists. It is well known a lot of results dealing with exact value of $h(G)$ for some special classes of graphs or upper and lower bounds on $h(G)$ for wider classes of graphs, for more see  \cite{Cat}, \cite{Hol}, \cite{CharWall}, \cite{Lai}, \cite{LiuXi}, \cite{Sar1}, \cite{Sar2}, \cite{Xio1}, \cite{Xio2}.

Surprisingly for the existence of hamiltonian paths in iterated line graphs there are not too much known results. In this paper we will consider hamiltonian paths in iterated line graphs, introduce the concept of hamiltonian path index and set the hamiltonian path index for some classes of graphs.

\section{ Terminology and Results}
 We say that $G$ is block-chain, if the block-cutvertex graph of $G$ is a path.

A 2-block is a 2-connected graph or a block of $G$ containing more than two vertices. An end block is a block of $G$ which contains at most one cutvertex, a block with at least two cutvertices of $G$ is non-end block. 

A dominating trail, a dominating closed trail in $G$ is a trail, a closed trail~$T$ in $G$ such that every edge of $G$ has at least one vertex on $T$, respectively. We say that an edge $e=uv$ of $G$ is  dominated by $u$ in $T$, if $u\in V(T)$ and $v\in V(G)\setminus V(T)$. 

Similarly as hamiltonian index  for hamiltonian cycles we define a hamiltonian path index for hamiltonian paths as follows. The hamiltonian path index, denoted by $h_p(G)$, is the minimum number $n$ such that $L^n(G)$ is traceable, it means that $L^n(G)$ contains a hamiltonian path. Clearly $h_p(G)\leq h(G)$ for all graphs except paths and $h_p(P)=0$ for every path $P$. Hence the hamiltonian path index of $G$ exists for any graph $G$. 

We are familiarized only with this arXiv paper \cite{Niu} by Niu et al. regar\-ding to the existence of a hamiltonian path in iterated line graphs. In this paper authors proved similar theorem (Theorem 4 in \cite{Niu}) as Xiong and Liu in \cite{Xio2}, that is  a characterization of G for which $L^n(G)$ has a hamiltonian path showing that there exists a nonempty set of subgraphs of $G$ satisfying 5~conditions, for $n\geq 2$. As applications, they use this characterization to give several upper bounds on the hamiltonian path index of a graph.

First we consider trees. Chartrand and Wall in \cite{CharWall} set the hamiltonian index of trees and we determine the hamiltonian path index for trees.

A branch in $G$ is a nontrivial path with end vertices of degree different from 2 and with internal vertices, if any, of degree 2. We denote by $B(G)$ the set of branches of $G$, by $CB(G)$ the subset of $B(G)$ in which every edge of every branch is a bridge of $G$, and by $CB_1(G)$ the subset of $CB(G)$ in which every branch has an end vertex of degree 1. Note that for a tree $T$ it holds that $B(T)$ is the same as $CB(T)$.

Xiong and Liu  in \cite{Xio2} define a parameter $k(G)$ for a graph $G$ depen\-ding on branches in $CB(G)$. Similarly we introduce a parameter $k(b)$ for every branch $b$ in $CB(G)$ as follows:

    $$k(b) =
    \begin{cases}
       |E(b)|+1, & \text{ if } b \in CB(G) \setminus CB_1(G),\\ 
       |E(b)|, & \text{ if } b \in CB_1(G).
    \end{cases} $$

\bigskip

\begin{figure}[h!]
    \centering
    \includegraphics[width = 8.5cm]{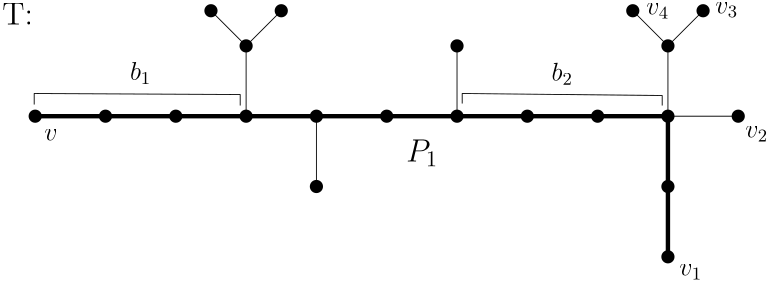}
    \caption{The $v_1v$-endpath $P_1$ (bold edges) in a tree $T$.}
    \label{koncova cesta P}
   \end{figure}

Let $T$ be a tree and not a path and $b_1$, $b_2$ be two arbitrary branches of $T$ such that $k(b_1)+k(b_2)$ is maximal. We denote by $\mathcal{P}$ a set of all unique $uv$-paths containing $b_1$, $b_2$, where $u, v$ are some two leaves in $T$.  For the illustration see Figure \ref{koncova cesta P} where $b_1$ ends in a leaf and $b_2$ not, hence every endpath $P_i$ in $\mathcal{P}$ ends in leaves $v$ and $v_i, i\in \{1,2,3,4\}$.

\bigskip

\noindent
For $P\in \mathcal{P}$ we denote by $B_{\mathit{P}}(T)$ all branches in $T$ contained in $P$. 

\newpage

Now we give a formula for the hamiltonian path index of trees.

\begin{theorem}
\label{Trees}
Let $T$ be a tree.  
\begin{itemize}
 \item $h_p(T)=0$ if and only if $T$ is a path;
 
\item $h_p(T)=1$ if and only if $T$ is a caterpillar and not a path;

 \item $T$ is not a caterpillar, then $$h_p(T)=\min _{P\in \mathcal{P} }\{\max\{\ k(b), \mbox{where } b\in CB(T)\setminus B_{\mathit{P}}(T) \}\}\geq 2.$$
\end{itemize}
\end{theorem}  

  As a generalization of trees we consider graphs with hamiltonian connected blocks. Clearly if $G$ has no 2-block, then $G$ is a tree and all blocks are bridges and hence hamiltonian connected. If $G$ contains a 2-block $B$, then we suppose that there is a hamiltonian path in $B$ between every two vertices of $B$. 

\begin{figure}[h!]
    \centering
    \includegraphics[width = 8.5cm]{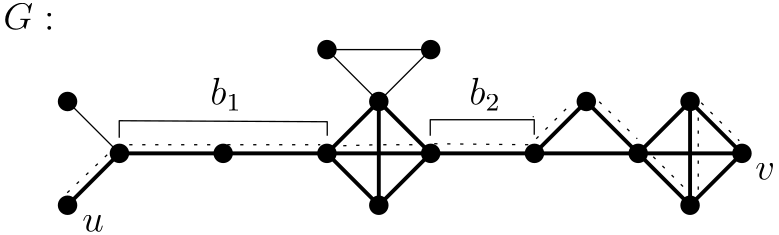}
    \caption{The $uv$-pseudopath (bold edges) from $uv$-path (dashed~edges) in~$G$.}
    \label{pseudocesta PP}
   \end{figure}

 Similarly as for trees, let $G$ be a graph with at least two branches in $CB(G)$ and $b_1, b_2$ be two arbitrary branches in $CB(G)$ such that $k(b_1)+k(b_2)$ is maximal.  Let $P$ be a $uv$-path in $G$ containing $b_1, b_2$, where $u, v$ are two vertices which are not cutvertices of $G$ and are contained in end blocks of~$G$. Then the $uv$-pseudopath $PP$ is a subgraph of $G$ which arises from some $uv$-path~$P$ by the replacement of all maximum paths (number of edges) in~$P$ containing in some 2-blocks of $G$,  by the whole 2-blocks.  Note that for given vertices $u,v$, a $uv$-path is not unique but corresponding $uv$-pseudopath is unique. For the illustration see Figure \ref{pseudocesta PP}.

\bigskip

\bigskip

\noindent
We denote by $\mathcal{PP}$ a set of all unique pseudopaths containing $b_1, b_2$ and for $PP\in \mathcal{PP}$ we denote by $B_{\mathit{PP}}(G)$ all branches in $CB(G)$ contained in $PP$. 

\newpage

We extend a formula for the hamiltonian path index of trees to graphs with hamiltonian connected blocks.

\begin{theorem}
\label{HamBlock}
Let $G$ be a graph with hamiltonian connected blocks and $M$ a~set of all leaves of $G$.

\begin{itemize}
 \item $h_p(G)=0$ if and only if $G$ is block-chain;
 \item $h_p(G)=1$ if and only if $G$ is not block-chain and all branches in $CB(G-M)$, if any, are contained on some path in $G-M$; 
 \item Otherwise $$h_p(G)=\min _{PP\in \mathcal{PP} }\{\max\{\ k(b), \mbox{where } b\in CB(G)\setminus B_{\mathit{PP}}(G) \}\}\geq 2.$$
\end{itemize}

\end{theorem}  

\bigskip

\section{Proofs}
Before proofs of main results of this paper we mention some needed known results.

%\begin{theorem}{(Chartrand and Wall, \cite{CharWall})} 
%\label{CharWall trees}
%Let $T$ be a tree. Then 
%$$h(T) = \max\{\ k(b), \mbox{where } b \in CB(T)\}.$$
%\end{theorem}

\begin{theorem}{(Harrary and Nash-Williams, \cite{HarNash})}
 \label{HarNash}
Let $G$ be a graph with at least three edges. Then $L(G)$ is hamiltonian if and only if $G$ has a~dominating closed trail.
\end{theorem}

\begin{theorem}{(Xiong and Zong, \cite{Xio3})} 
\label{Xio3}
Let $G$ be a graph with at least one edge. Then $L(G)$ is traceable if and only if $G$ has a dominating trail.
\end{theorem}

\medskip

\noindent
PROOF OF THEOREM \ref{Trees}

\begin{proof}
  Clearly, $h_p(T)=0$ if and only if $T$ is a path. 

\medskip

Suppose that $T$ is not a path. Let  $b_1, b_2$ be branches in $T$ such that $k(b_1)+k(b_2)$ is maximal. Note that it is possible to choose arbitrary pair of branches $b_i, b_j$ for which $k(b_i)+k(b_j)$ is maximal. 

 Let $P \in \mathcal{P}$ be an endpath for which maximal $k(b)$ is minimal for all $b\notin B_{\mathit{P}}(T)$. Basically,  we prove that the hamiltonian path index of a tree $T$ is exactly this maximum value of the parameter $k(b)$ for  all branches which are not in $B_{\mathit{P}}(T)$ (i.e. $b$ is not contained in $P$).  For convenience, we set $m=\min _{P\in \mathcal{P} }\{\max\{\ k(b), \mbox{where } b\in CB(T)\setminus B_{\mathit{P}}(T) \}\}$. 

If $T$ is a caterpillar then $P$ is the dominating trail in $T$ and $L(T)$ has a~hamiltonian path by Theorem \ref{Xio3}. Hence $h_p(T) =1$.  On the other side, if $h_p(T)=1$, then $T$ has a dominating trail  by Theorem \ref{Xio3} and $T$ has to be a~caterpillar.
  
\medskip

Now suppose that $T$ is not caterpillar. Hence $m\geq 2$.

\smallskip

\noindent
 First we show that $h_p(T) \leq m$. We proceed similarly as in \cite{CharWall} and start with the following observations which follow immediately.

\bigskip

\emph{Observation 1:}
\begin{mathitem}
\item[]
From the manner in which the line graph is defined, each 2-block of $L(T)$ is produced from a subgraph of $T$ which is a star. Hence each 2-block of $L(T)$ is complete and therefore hamiltonian and also even hamiltonian connected. 

In $L^2(T)$ the 2-blocks are produced using those maximal unions of hamiltonian cycles in 2-blocks of $L(T)$ which create a closed trail in~$L(T)$. Thus each 2-block of $L^2(T)$ is hamiltonian by Theorem~\ref{HarNash}. 

Unfortunately the 2-blocks of $L^2(T)$ are not hamiltonian connected in general. But in $L^2(T)$ the 2-blocks are also produced using those maximal unions of hamiltonian cycles and two hamiltonian paths in 2-blocks of $L(T)$ which create a trail $DT$ in $L(T)$.  Let $u, v$ be vertices of odd degree in $DT$. We denote by $D_u$, $D_v$ a set of dominated edges by $u$, $v$, respectively.

Then there exists a hamiltonian path in $L^2(T)$ between a vertex, which corresponds either to some edge of $D_u$ if $D_u\neq\emptyset$ or to some edge $uu'\in E(DT)$ if $D_u=\emptyset$ , and a vertex, which corresponds to some edge of $D_v$ if $D_v\neq\emptyset$ or to some edge $vv'\in E(DT)$ if $D_v=\emptyset$, by Theorem~\ref{Xio3}. 

The same statements for hamiltonian cycles and paths can now be made about the graphs $L^n(T), n\geq 3$. 
\end{mathitem}

\bigskip

\emph{Observation 2:} 
 \begin{mathitem}
\item[]
Let $b$ be a branch in $CB_1(T)$.  

If $k(b)\leq m$, then there is no branch corresponding to $b$ in $L^m(T)$ and there is a 2-block in $L^m(T)$ which is produced from line graph operation on a subgraph of $T$ containing~$b$. 

If $k(b)>m$, then there is a branch in $CB_1(L^m(T))$ corresponding to $b$ in $L^m(T)$ having length $k(b)-m$.  
 \end{mathitem}

\bigskip

\emph{Observation 3:}
\begin{mathitem}
\item[]
Let $b$ be a branch in $CB(T)\setminus CB_1(T)$. 

If $k(b)\leq m$, then there is neither a branch nor a cutvertex corresponding to $b$ and there is a 2-block in $L^m(T)$ which is produced from line graph operation on a subgraph of $T$ containing $b$.

 If $k(b)=m+1$, then there is a cutvertex corresponding to $b$ in $L^m(T)$ which is produced from line graph operation on $b$. 

If $k(b)>m+1$, then there is a branch in $CB(L^m(T))\setminus CB_1(L^m(T))$ corresponding to $b$ in $L^m(T)$ having length $k(b)-m-1$.  
\end{mathitem}

\bigskip

With respect to the choice of $P$ and Observation 2 and 3 it holds that $L^m(T)$ is block-chain. If $L^m(T)$ is a 2-block, then $L^m(T)$ is hamiltonian by Observation 1. Hence it contains a hamiltonian path and $h_p(T) \leq m$.

Now suppose that $L^m(T)$ is not a 2-block. Let $v_1, v_2$ be both end vertices of $P$. For $i=1,2$, if $v_i$ is contained in a~branch $b^i$ in $CB_1(T)$ with $k(b^i)> m$, then we set the path $P^i=L^m(b^i)$; if $v_i$ is contained in a~branch $b^i$ in $CB_1(T)$ with $k(b^i)\leq m$, then we denote by $B^i$ the 2-block in $L^m(T)$ corresponding to $b^i$  which is end block and hamiltonian (see Observation 1 and 2) and we set $P^i$ a hamiltonian path in $B^i$ which arises from a hamiltonian cycle in $B^i$ deleting an edge containing the unique cutvertex of $L^m(T)$ in $B^i$.  

We denote by $\mbox{br}_i\in B_P(T)\setminus CB_1(T)$ with $k(\mbox{br}_i)> m+1$, $i=1,2,...k$. We set the paths $P_i=L^m(\mbox{br}_i)$ which are also branches in $L^m(T)$ (see Observation~3). 

Note that $b_1$, $b_2$ could be included in $b^1, b^2$ or $\mbox{br}_i$, $i=1,2,...k$, if $k(b_1)$, $k(b_2)$ is sufficiently large (see Observation 2 and 3).  

Now consider all non-end 2-blocks $B_i$ in $L^m(T)$, $i=1,2,...,l$. Every $B_i$ contains exactly 2 cutvertices $c_i,c'_i$ of $L^m(T)$(because of $L^m(T)$ is block-chain) and $B_i$ contains a hamiltonian path $P_{B_i}$ with $c_i, c'_i$ as end vertices of $P_{B_i}$ (see Observation 1).  Note that $c_i$, $c'_i$ corresponds to some edge of $D_u$, $D_v$ in $L^{m-1}(T)$ by Observation 1, respectively. 

\begin{figure}[h!]
    \centering
    \begin{subfigure}{0.4\textwidth}
        \includegraphics[width=\linewidth]{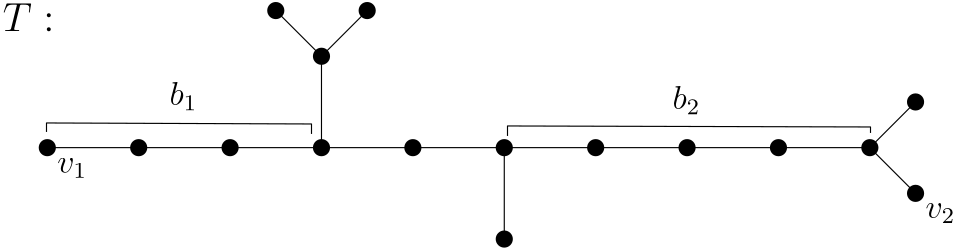}
        \label{figure_2a}
    \end{subfigure}
    \hspace*{2cm} % mezera mezi obrázky
    \begin{subfigure}{0.4\textwidth}
        \includegraphics[width=\linewidth]{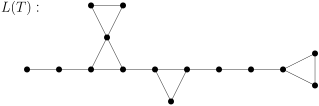}
        \label{figure_2b}
    \end{subfigure}
    \begin{subfigure}{0.5\textwidth}
        \includegraphics[width=\linewidth]{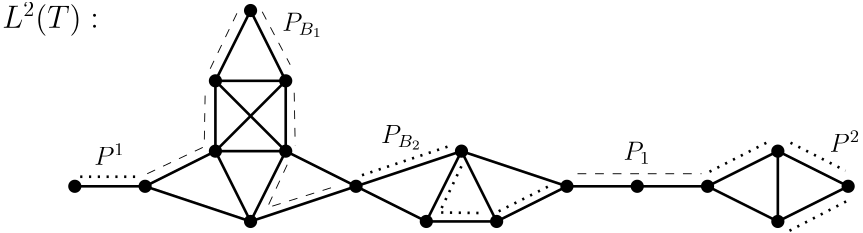}
        \label{figure_2c}
    \end{subfigure}
    \caption{~}
    \label{figure_2}
\end{figure}

The constructions from the last four paragraphs above are illustrated on Figure \ref{figure_2} for a tree $T$, where $m=2$. Resulting hamiltonian path in $L^2(T)$ is marked with dots and dashes.

Finally $$P_h=P^1\cup P^2\cup_{i=1}^k P_i \cup_{j=1}^l P_{B_i}$$ is a hamiltonian path in $L^m(T)$. Hence $L^m(T)$ is traceable and $h_p(T) \leq m$.

\bigskip

\noindent
Now we show that $h_p(T) \geq m$. 

\smallskip

Let $b^*, b^{**} \in B_P(T)$ be two branches such that $k(b^*)\geq m$, $k(b^{**})\geq m$, and  $b' \in CB(T)\setminus B_P(T)$ be a branch such that $k(b')=m$. Let $P^1=a_1 a_2... a_k$ be a unique path in $T$ where $a_1\in V(b^*), a_k\in V(b^{**})$ and $a_i\notin V(b^*)\cup V(b^{**})$ for $i=2,3,...,k-1$, and $P^2=b_1 b_2... b_l$ be a unique path in $T$ where $b_1\in V(b'), b_l \in V(P)$ and $b_i\notin V(b')\cup V(P)$ for $i=2,3,...,l-1$. Note that $P^1$, $P^2$ could be trivial, i.e. $P^1=a_1=a_k$, $P^2=b_1=b_l$, respectively. For illustration see Figure \ref{figure_3}a).

\bigskip

\begin{figure}[h!]
    \centering
    \begin{subfigure}{0.45\textwidth}
        \includegraphics[width=\linewidth]{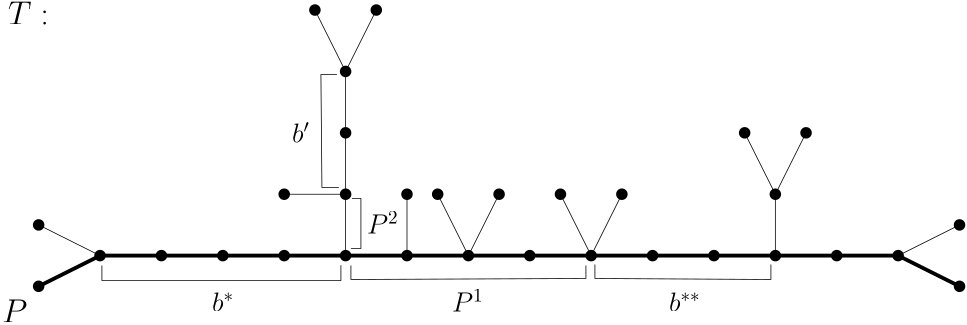}
        \caption*{\scriptsize a) Branches $b^*, b^{**}, b'$ and paths $P^1, P^2$ in $T$.}
        \label{figure_3a}
    \end{subfigure}
    \hspace*{2cm} % mezera mezi obrázky
    \begin{subfigure}{0.3\textwidth}
        \includegraphics[width=\linewidth]{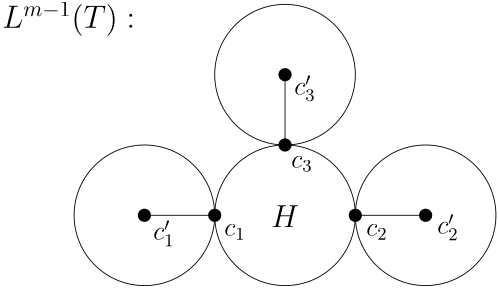}
        \caption*{\scriptsize b) The structure of $L^{m-1}(T).$}
        \label{figure_3b}
    \end{subfigure}
    \caption{~}
    \label{figure_3}
\end{figure}

\bigskip

Note that the path $P$ is chosen such that we want to minimize a maximum value of $k(b)$ for branches not in $P$. It means that $P$ could contain instead of $b_1, b_2$ another branches with $k(b)\geq m$. With respect to the previous fact  branches $b^*, b^{**},b'$ exist such that $V(P^1)\cap V(P^2)\neq \emptyset$. Otherwise we find another path $P'$ containing all branches with  $k(b)\geq m$ as $P$ and also $b'$, a contradiction with the fact that we minimize a maximum value of $k(b)$ for branches not in $P$.

Since $T$ is not a caterpillar, then $k(b^*)\geq m>1$, $k(b^{**})\geq m>1$ and $k(b')=m>1$, and hence $L^{m-1}(T)$ contains at least three cutvertices $c_1, c_2, c_3$ which arise from $b^*, b^{**}, b'$ by successive repeating line graph operation on $T$. Moreover since $V(P^1)\cap V(P^2)\neq \emptyset$, there exists a subgraph $H$ of $L^{m-1}(T)$ such that $c_i\in V(H)$ and there exists $c'_i\in N_{L^{m-1}(T)}(c_i)\setminus V(H)$ for $i=1, 2, 3$. Obviously there is not a path in $L^{m-1}(T)$ containing $c'_1, c'_2, c'_3$ (see Figure~\ref{figure_3}b)). Hence $h_p(T) \geq m$.
\end{proof}

\medskip

\noindent
PROOF OF THEOREM \ref{HamBlock}

\begin{proof}
First suppose that $G$ is block-chain. Then clearly there is a hamiltonian path in $G$, hence $h_p(G) = 0$.  We use a hamiltonian path in every non-end block between its 2 cutvertices of $G$ and a hamiltonian path in eve\-ry end-block between its cutvertex of $G$ and its arbitrary other vertex.  If $G$ is not block-chain, then clearly there is not a hamiltonian path in $G$, hence $h_p(G)>0$. The first part of Theorem \ref{HamBlock} is proved.

Now  $G$ is not block-chain.

Suppose that all branches in $CB(G-M)$, if any, are contained on some path in $G-M$. Let $P$ be a $uv$-path containing all branches in $CB(G-M)$ where $u, v$ are in different end-blocks of $G-M$ and are not cutvertices of $G-M$. Let $PP$ be a corresponding pseudopath from $P$. Then $PP$ is block-chain and contains a hamiltonian path (see the paragraph above). This hamiltonian path in $PP$ and hamiltonian cycles in all other blocks of $G-M$ not in $PP$ (note that these blocks are 2-blocks) are a dominating trail in $G$ and $L(G)$ has a hamiltonian path by Theorem \ref{Xio3}, hence $h_p(G) = 1$. If $G-M$ contains at least three branches in $CB(G-M)$ not in a path, then clearly $G$ has not a dominating trail and by Theorem \ref{Xio3} there is not a hamiltonian path in $L(G)$, hence  $h_p(G)>1$. The second part of Theorem \ref{HamBlock} is proved.

Now suppose that $G-M$ contains at least three branches in $CB(G-M)$ not in a path and $h_p(G)\geq 2$. Hence there is at least one branch in $CB(G)\setminus B_{\mathit{PP}}(G)$ for every pseudopath $PP$. We continue on the same principle as in the proof of Theorem \ref{Trees}. Basically, each 2-block of the line graph of a tree is hamiltonian connected because it is complete (see Observation~1 in the proof of Theorem  \ref{Trees}) and each 2-block of a graph $G$ is hamiltonian connected from assumptions of Theorem \ref{HamBlock}, which is a crucial point in the proof of Theorem~\ref{Trees}. We replace paths in proof of Theorem \ref{Trees} with corresponding pseudopaths and continue with the same arguments as in the proof of Theorem \ref{Trees}.
\end{proof}

\section{Conclusion} 
We set the exact value of the hamiltonian path index of trees and graphs with hamiltonian connected blocks. Sara\v zin in \cite{Sar1} proved that the hamiltonian index is the same for trees and graphs with hamiltonian 2-blocks. Maybe surprisingly this is not true for the hamiltonian path index.  

\bigskip

    \begin{figure}[h!]
    \centering
    \includegraphics[width = 8.5cm]{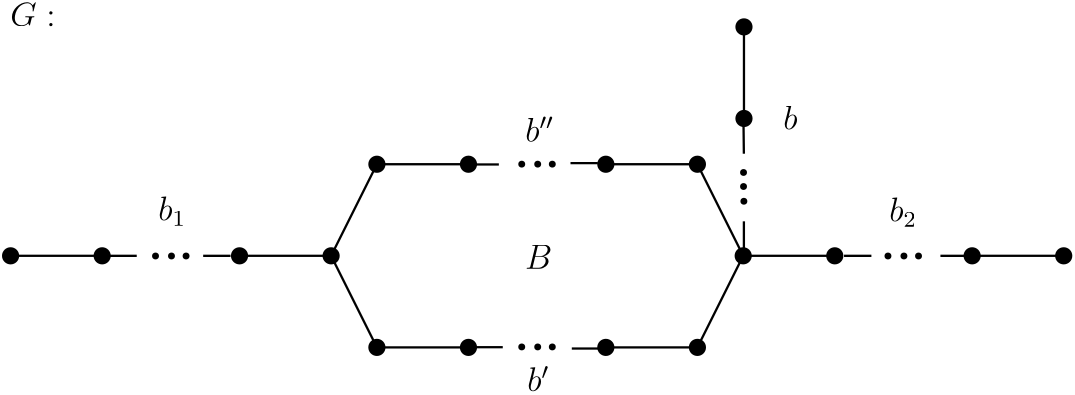}
    \caption{The example of graphs $G$ with hamiltonian 2-blocks such that $k(b)=|E(b)|=m$ and $h_p(G)>m$ for sufficiently large $b_1, b_2, b', b''$.}
    \label{Counterexample}
   \end{figure}

\bigskip

Consider the graph $G$ on Figure \ref{Counterexample} for fixed $m$. If we use Theorem~\ref{HamBlock} for $G$, then the hamiltonian path index of $G$ should be $k(b)=m$ but really $h_p(G)$ depends on branches $b', b''$ of the hamiltonian block $B$ and branches $b_1, b_2$. If the lenghts of branches $b_1, b_2, b', b''$ are sufficiently large, then also $h_p(G)$ is sufficiently large not depending on $m$. Note that all these branches are on a~pseudopath $PP$ for which the minimum in the formula in Theorem~\ref{HamBlock} is reached for $m\geq2$.   

\bigskip

On the other side the condition on hamiltonian connectedness of all blocks is too strong. Following proofs of Theorem \ref{Trees} and Theorem \ref{HamBlock} we can state this theorem as a corollary.

\begin{theorem}
 \label{HamBlock2}
Let $G$ be a graph and $M$a set of all leaves of $G$. 

\medskip

\noindent
 a) Let $G$ be block-chain. 

\noindent
If all non-end blocks of $G$ have a hamiltonian path between their cutvertices and both end blocks of $G$ have a hamiltonian path between their cutvertex and other their vertex, then $h_p(G)=0$.

\medskip

\noindent
 b)  Let $G$ be not block-chain. 
Let $PP$ be a pseudopath in $G$ either such that contains all branches in $CB(G-M)$ and $m=1$, or for which  $m=\min _{PP\in \mathcal{PP} }\{\max\{\ k(b), \mbox{where } b\in CB(G)\setminus B_{\mathit{PP}}(G) \}\}$ is reached.

\noindent
 If all blocks not in $PP$ are hamiltonian, all non-end blocks of $PP$ have a~hamiltonian path between their cutvertices and both end blocks of $PP$ have a hamiltonian path between their cutvertex and other their vertex, then $h_p(G)=m$.
 \end{theorem}

\bigskip

We think in a connection with graphs $G$ on Figure \ref{Counterexample} that the~finding of a~hamiltonian path index of not only graphs with hamiltonian 2-blocks leads to branch bonds (for more see \cite{XioBro}) where also even branch bonds will play important role (not only odd branch bonds as for the hamiltonian index). We leave it as an open problem.

\bigskip

\noindent {\bf Acknowledgements}. 
The second author was partially supported by the internal grant SGS-2025-007 of the University of West Bohemia in Pilsen.

\end{document}